\documentclass[12pt,leqno]{article}
\usepackage{amsmath,amssymb,amscd,latexsym}

\newcommand{\A}{{\mathbb{A}}}

\newcommand{\F}{{\mathbb{F}}}
\newcommand{\Lee}{\mathbb{L}}
\newcommand{\Pa}{{\mathbb{P}}}

\newcommand{\Z}{{\mathbb{Z}}}

\newcommand{\mot}{\mathrm{mot}}
\newcommand{\sep}{\mathrm{sep}}

\newcommand{\Spec}{\mathrm{Spec}\,}
\newcommand{\Pic}{\mathrm{Pic}}

\newcommand{\Var}{\mathrm{Var}}
\newcommand{\Dh}{{\mathcal D}}

\newcommand{\Lh}{{\mathcal L}}

\newcommand{\Sh}{{\mathcal S}}
\newcommand{\betr}[1]{|#1|}
\newcommand{\ohne}{\setminus}
\newcommand{\dis}{\displaystyle}

\newtheorem{theorem}{Theorem}

\newenvironment{proof}{\noindent {\bf Proof}}{\mbox{}\hfill$\Box$}
\begin{document}
\title{ A motivic version of Pellikaan's two variable zeta function}
\author{Francesco Baldassarri \and Christopher Deninger \and Niko Naumann}
\date{\ }
\maketitle

Consider the Grothendieck ring $K_0 (\Var_K)$ of varieties over a field $K$. It is the abelian group generated by the isomorphism classes $[X]$ of algebraic $K$-schemes subject to the relations
\[
[X] = [Y] + [X \ohne Y]
\]
for any closed subscheme $Y$ of $X$. The product is defined by the formula
\[
[X_1] \cdot [X_2] = [X_1 \times_K X_2] \; .
\]
The class of $\A^1$ plays a special role. It is denoted by $\Lee = [\A^1]$.

For a variety $X / K$, Kapranov \cite{K} defined the motivic zeta function of $X$ as the power series
\[
Z_{X,\mot} (T) = \sum^{\infty}_{n=0} \; [X^{(n)}] T^n \quad \mbox{in} \; K_0 (\Var_K) [ \betr{T}] \; .
\]
Here $X^{(n)} = X^n / \Sh_n$ is the $n$-fold symmetric product, $X^{(0)} = \Spec(K)$. Moreover, for a homomorphism of commutative rings with unity,
\[
\mu : K_0 (\Var_K) \longrightarrow A,
\]
a so called motivic measure, he considers the series
\[
Z_{X,\mu} (T) = \sum^{\infty}_{n=0} \mu ([X^{(n)}]) T^n \quad \mbox{in} \; A [\betr{T}] \; .
\]
The rationality of these power series is an interesting issue. See \cite{DL} and \cite{LL}.
For $X/\F_q$ and $\mu:K_0(\Var_{\F_q})\rightarrow\Z$ defined by $\mu([Z])=|Z(\F_q)|$, 
$Z_{X,\mu}(T)$ equals the usual zeta function of $X/\F_q$, i.e.:

\[
Z_X(T):= \mbox{exp}(\sum_{\nu\ge1} |X(\F_{q^\nu})|\frac{t^\nu}{\nu})=\sum_{\nu\ge 0} |X^{(\nu)}(\F_q)|t^{\nu}.
\]

This is the special case of constant coefficients concentrated in degree zero of
\cite{D}, lemme 4.11. It can also be proved by a direct combinatorial argument without
using any cohomological expression for $Z_X(T)$.

Concerning curves, Kapranov proves the following result in \cite{K}:

\begin{theorem}
  [Kapranov] Let $X$ be a smooth projective geometrically connected curve of genus $g$ over $K$ with a degree $1$ line bundle. Assume that $A$ is a field and that $\Lee_{\mu} := \mu (\Lee) \neq 0$ in $A$. Then we have
\[
Z_{X,\mu} (T) = \frac{P_{X,\mu} (T)}{(1-T) (1 -\Lee_{\mu} T)}
\]
for a polynomial $P_{X,\mu} (T)$ of degree $2g$ in $A [T]$. Moreover the following functional equation holds:
\[
Z_{X,\mu} (T) = \Lee^{g-1}_{\mu} T^{2g-2} Z_{X,\mu} (\Lee^{-1}_{\mu} T^{-1} ) \; .
\]
\end{theorem}

For $\mu : K_0 (\Var_{\F_q}) \to \Z$ as above the theorem asserts the well known rationality and symmetry properties of the usual zeta function of a curve over a finite field.

We quickly recall the construction, due to R. Pellikaan, of a two-variable zeta function
(see \cite{P} for more details):

Let $X/\F_q$ be a curve. Its (usual) zeta function can be written as

\[
Z_X(T)=\sum_{D\ge 0}T^{\deg(D)},
\]

the sum being extended over all effective divisors on $X$. 
Summing over divisor-classes $\Dh=[D]$ this becomes

\[
Z_X(T)=\sum_{\Dh}\frac{q^{h^0(\Dh)}-1}{q-1}T^{\deg(\Dh)}.
\]

Here we write $h^0(\Dh)=$ dim $H^0(X,\mathcal{O}(D))$ for any divisor $D$ with $\Dh=[D]$.
The two-variable zeta function is obtained by substituting a variable $u$ for $q$ in
this expression:

\[
Z_X(T,u)=\sum_{\Dh}\frac{u^{h^0(\Dh)}-1}{u-1}T^{\deg(\Dh)}.
\]

It is a power-series in $u$ and $T$ with integer coefficients. Pellikaan then proves, among other things, rationality:

\[
Z_X(T,u)=\frac{P_X(T,u)}{(1-T)(1-uT)},
\]

for a suitable $P_X(T,u)\in\Z[T,u]$. It was suggested by a question of J. Lagarias and
E. Rains and proven in \cite{N} that $P_X(T,u)\in\Z[T,u]$ is absolutely irreducible.

A two-variable zeta function may be defined in the motivic case as well and as we shall see, it has similar properties as $Z_X (T,u)$. For a curve $X / K$ as in Kapranov's theorem let $\Pic^n_{X/K}$ be the Picard variety of degree $n$ line bundles on $X$. Let $\Pic^n_{\ge \nu}$ be the closed subvariety in $\Pic^n_{X/K}$ of line bundles $\Lh$ with $h^0 (\Lh) \ge \nu$. The algebraic $K$-scheme $\Pic^n_{\nu} = \Pic^n_{\ge \nu} \ohne \Pic^n_{\ge \nu + 1}$ defines a class in $K_0 (\Var_K)$
\[
[\Pic^n_{\nu}] = [\Pic^n_{\ge \nu}] - [\Pic^n_{\ge \nu+1}] \; .
\]
The two-variable motivic zeta function of the curve $X / K$ is defined as the formal power series
\[
Z_{X,\mot} (T,u) = \sum_{n , \nu \ge 0} [\Pic^n_{\nu}] \frac{u^{\nu}-1}{u-1} T^n \quad \mbox{in} \; K_0 (\Var_K) [|u,T|] \; .
\]
For a motivic measure $\mu$ with values in $A$ we set:
\[
Z_{X,\mu} (T,u) = \sum_{\nu,\mu \ge 0} \mu [\Pic^n_{\nu}] \frac{u^{\nu}-1}{u-1} T^n \quad \mbox{in} \; A [|u,T|] \; .
\]

Here, again, for $X/\F_q$ and $\mu:K_0(\Var_{\F_q})\rightarrow \Z$ as above we have
$Z_{X,\mu}(T,u)=Z_X(T,u)$, because for any $n, \nu\ge 0$ the set of $\F_q$-rational 
divisor classes $\Dh$ with $\deg(\Dh)=n$ and $h^0(\Dh)=\nu$ are in bijection with
the $\F_q$-rational points of $\Pic^n_{\nu}$. This follows from a general result on
the relative Picard functor, \cite{BLR}, 8.1, Prop. 4, and the fact that the 
Brauer group of a finite field is trivial.

We note that $Z_{X,\mu}(T,u)$ is constructed very much the same way as Pellikaan did:
The natural morphism 

\[
X^{(n)}\longrightarrow \Pic^n_{X/K}
\]

is such that $X^{(n)}\times_{\Pic^n_{X/K}} \Pic^n_{\nu}\simeq \Pa^{\nu-1}\times \Pic^n_{\nu}$
which implies that (c.f. \cite{K}, Prop. 1.2.3)

\[
[X^{(n)}]=\sum_{\nu\ge 0} [\Pic^n_{\nu}]\frac{\Lee^{\nu}-1}{\Lee-1}\quad\mbox  { in } K_0(\Var_K)
\]
and hence

\[
Z_{X,\mot}(T)=\sum_{n,\nu\ge 0} [\Pic^n_{\nu}] \frac{\Lee^{\nu}-1}{\Lee-1}T^n.
\]

Just as Pellikaan substituted a variable for the integer $q$, we substitute a
variable for the class $\Lee$ in $K_0(\Var_K)$.

We explain a convenient way of writing $Z_{X,\mu}(T,u)$ using motivic integration:
Let $K^{\sep}$ denote a fixed separable closure of $K$ and let $Z$ be an algebraic scheme over $K$. A constructible $A$-valued function on $Z$ in the sense of \cite{K} (1.2) is a function
\[
f : Z (K^{\sep}) \longrightarrow A
\]
which can be written in the form 
\[
f = \sum^n_{i=1} a_i \chi_{W_i (K^{\sep})}
\]
where $a_i \in A$ and the $W_i \subset Z$ are closed subschemes. Here $\chi$ denotes the characteristic function of a set.

The integral of $f$ with respect to $\mu$ is defined to be
\[
\int_Z f \, d\mu = \sum^n_{i=1} a_i \mu (W_i) \; .
\]
In our case, the function
\[
\Pic^n_{X/K} (K^{\sep}) \longrightarrow A [\betr{u}] \; , \; \Dh \longmapsto \frac{u^{h^0 (\Dh)}-1}{u-1}
\]
defines a constructible $A[\betr{u}]$-valued function on $\Pic^n_{X/K}$ and we may write $Z_{X , \mu} (T,u)$ in the form:
\begin{equation}
  \label{eq:1}
  Z_{X,\mu} (T,u) = \sum_n \left( \int_{\Pic^n_{X/K}} \frac{u^{h^0 (\Dh)}-1}{u-1} \, d \mu (\Dh) \right) T^n \; .
\end{equation}
Here $\mu$ is viewed as an $A [\betr{u}]$-valued measure using the inclusion $A \hookrightarrow A [\betr{u}]$. 

Introducing a suitable notion of convergent integrals over $K$-schemes which are not of finite type we could also write
\[
Z_{X,\mu} (T,u) = \int_{\Pic_{X/K}} \frac{u^{h^0 (\Dh)}-1}{u-1} T^{\deg \Dh} \, d\mu (\Dh)
\]
where now $\mu$ is viewed as an $A [\betr{u}] [\betr{T}]$-valued measure. However we will do with (\ref{eq:1}) in the sequel. The following result is proved in the same way as Pellikaan's original theorem for $Z_X(T,u)$. 

\begin{theorem}\label{functionalequation}
  Let $K$ be a field and $X / K$ a smooth projective geometrically irreducible curve of genus $g$ which admits a line bundle of degree one. Then we have\\
a) $Z_{X,\mu} (T,u) = \frac{P_{X,\mu} (T,u)}{(1-T) (1-uT)} \quad \mbox{in} \; A [T, u, (1-T)^{-1} , (1- uT)^{-1}]$\\
where $P_{X,\mu} (T,u) \in A [T,u]$. \\
b) $\dis P_{X,\mu} (T,u) = \sum^{2g}_{i=0} P_i (u) T^i \quad \mbox{with} \; P_i(u) \in A [u]$, where \\
$P_0 (u) = 1 , P_{2g} (u) = u^g , \deg P_i (u) \le 1 + \frac{i}{2}$.\\
c) $\dis Z_{X,\mu} (T, u) = u^{g-1} T^{2g-2} Z_{X,\mu} \left( \frac{1}{Tu} , u \right)$ i.e. $P_{2g -i} (u) = u^{g-i} P_i (u)$.\\
d) $P_{X,\mu}(1,u)=\mu(\Pic^0_{X/K})\in A\subset A[u]$.
\end{theorem}

\begin{proof}
  For $g = 0$ we have $Z_{X,\mu} (T,u) = (1-T)^{-1} (1 - uT)^{-1}$, so the assertions are clear. For $g \ge 1$ we have:
  \begin{eqnarray*}
    Z_{X,\mu} (T,u) & = & \sum_{0 \le n \le 2g -2} T^n \int_{\Pic^n_{X / K}} \frac{u^{h^0 (\Dh)} - 1}{u-1} d\mu (\Dh) \\
 && + \sum_{n > 2g-2} \frac{u^{n+1-g}-1}{u-1} T^n \int_{\Pic^n_{X/K}} \, d\mu (\Dh) \\
& = & \sum_{0 \le n \le 2g-2} T^n \int_{\Pic^n_{X/K}} \frac{u^{h^0 (\Dh)}}{u-1} \, d\mu (\Dh) \\
 && + \sum_{n > 2g-2} \frac{u^{n+1-g}}{u-1} \mu (\Pic^n_{X/K}) T^n - \sum_{n \ge 0} \frac{\mu (\Pic^n_{X/K})}{u-1} T^n \; .
  \end{eqnarray*}
Since there exists a degree one line bundle on $X$, we have $\Pic^n_{X /K} \cong \Pic^0_{X / K}$. Hence we get
\begin{eqnarray*}
  Z_{X,\mu} (T,u) & = & \frac{1}{u-1} \sum_{0 \le n \le 2g-2} T^n \int_{\Pic^n_{X/K}} u^{h^0 (\Dh)} \, d\mu (\Dh) \\
 && + \frac{\mu (\Pic^0_{X/K})}{u-1} \left( u^{1-g} \frac{(uT)^{2g-1}}{1-uT} - \frac{1}{1-T} \right) \\
& = & \frac{1}{u-1} \int_{\Pic^{[0,2g-2]}_{X/K}} u^{h^0 (\Dh)} T^{\deg \Dh} \, d\mu (\Dh) \\
 & & + \frac{\mu (\Pic^0_{X/K})}{u-1} \left( \frac{u^g T^{2g-1}}{1-uT} - \frac{1}{1-T} \right) \; .
\end{eqnarray*}
In the last formula we have set $\Pic^{[0,2g-2]}_{X/K} = \coprod_{0 \le n \le 2g-2} \Pic^n_{X/K}$.
The rest of the proof is as in \cite{P} but for convenience we give the details:
We set

\begin{eqnarray*}
G(T,u)=\int_{\Pic^{[0,2g-2]}_{X/K}} u^{h^0 (\Dh)} T^{\deg \Dh} \, d\mu (\Dh) & & \\
F(T,u)=\mu (\Pic^0_{X/K}) \left( \frac{u^g T^{2g-1}}{1-uT} - \frac{1}{1-T} \right). & & 
\end{eqnarray*}

So we have $(u-1)Z_{X,\mu}(T,u)=G(T,u)+F(T,u)$. A direct computation shows that $F(T,u)$ satisfies
the functional equation. To check the functional equation for $G(T,u)$ we observe that sending 
a line bundle $\Lh$ to $\Lh'=\Lh^{-1}\otimes\omega_{X/K}$, where $\omega_{X/K}$ is the canonical bundle, defines an
involution on the relative Picard functor of $X/K$ because formation of the canonical bundle 
commutes with arbitrary base-change.
The theorem of Riemann-Roch together with Serre-duality on our curve implies that the corresponding involution on the $K$-scheme $\Pic_{X/K}=\amalg_{n\in\Z}\Pic^n_{X/K}$ induces isomorphisms for any $n, \nu\ge 0$:

\[
\Pic^n_{\nu}\stackrel{\simeq}{\longrightarrow} \Pic^{2g-2-n}_{\nu-n+g-1}\; , \;  \Dh\mapsto\Dh'
\]

Using this we compute:
\begin{eqnarray*}
& & u^{g-1}T^{2g-2}G(\frac{1}{Tu},u)\\
 & = & \int_{\Pic^{[0,2g-2]}_{X/K}} u^{h^0(\Dh)+g-1-\deg(\Dh)}T^{2g-2-\deg(\Dh)} \, d\mu (\Dh)\\
& = & \int_{\Pic^{[0,2g-2]}_{X/K}} u^{h^0(\Dh')}T^{\deg(\Dh')} \, d\mu (\Dh')\\
& = & G(T,u). 
\end{eqnarray*}

So the functional equation holds for $F(T,u)$ and $G(T,u)$, hence for $Z_{X,\mu}(T,u)$, i.e.
we have proved the first assertion of part c). Consider
\begin{eqnarray} \label{eq1}
 &   & Q(T,u):=(1-uT)(1-T)(F(T,u)+G(T,u))\\
 & = & \mu(\Pic^0_{X/K})((1-T)u^gT^{2g-1}-(1-uT))+(1-uT)(1-T)G(T,u).\nonumber
\end{eqnarray}

This is a polynomial in $u$ and $T$ with coefficients in $A$.
As

\begin{eqnarray*}
 &   & G(T,1)= \int_{\Pic^{[0,2g-2]}_{X/K}} T^{\deg(\Dh)} \, d\mu (\Dh)\\
 &  = & \sum_{n=0}^{2g-2} \mu(\Pic^n_{X/K}) T^n  = \mu(\Pic^0_{X/K})\frac{T^{2g-1}-1}{T-1}
\end{eqnarray*}

one computes $Q(T,1)=0$, hence $Q(T,u)=(u-1)P_{X,\mu}(T,u)$ for some $P_{X,\mu}(T,u)\in A[T,u]$.
Since $u-1$ is not a zero-divisor in $A[T,u]$ we find

\[
P_{X,\mu}(T,u)=(1-uT)(1-T)Z_{X,\mu}(T,u),
\]

i.e. we have proved part a). Clearly, $G(T,u)$ has degree at most
$2g-2$ in the variable $T$. So from (\ref{eq1}) the same holds true for $P_{X,\mu}(T,u)$.
We also have

\[
Q(1,u)=\mu(\Pic^0_{X/K})(u-1)=(u-1)P_{X,\mu}(1,u),
\]

which proves part d). It is left to the reader to check that the functional equation
is indeed equivalent to $P_{2g-i}(u)=u^{g-i}P_i(u)$ for all $i\ge 0$.
To prove part b) we have

\[
Q(0,u)=-\mu(\Pic^0_{X/K})+G(0,u).
\]

Observe that $\Pic^0_{\ge 1}=\Pic^0_1$ is the zero section of the abelian variety 
$\Pic^0_{X/K}$ and hence $\mu(\Pic^0_1)=1=\mu(\Pic^0_{X/K})-\mu(\Pic^0_0)$ from which we get

\begin{eqnarray*}
 &  &  Q(0,u)=-\mu(\Pic^0_{X/K})+\int_{\Pic^0_{X/K}} u^{h^0(\Dh)} \; d\mu(\Dh) \\
 &  = & -\mu(\Pic^0_{X/K})+\mu(\Pic^0_0)+u = u-1 \;,
\end{eqnarray*}

i.e. $P_0(u)=1$ and from the functional equation $P_{2g}(u)=u^g$ also.
The last assertion, $\deg\; P_i(u)\leq 1+i/2$, follows from Clifford's theorem
(c.f. \cite{H}, IV, thm. 5.4), which can be formulated as asserting that if a
divisor class $\Dh$ has $h^0(\Dh)\geq \mbox{max}\{1,\deg(\Dh)-g+2\}$ then $h^0(\Dh)\leq\frac{1}{2}\deg(\Dh)+1$.
This implies that for any $n, \nu\ge 0$, if $\nu\ge \mbox{max}\{1,n-g+2\}$ and $\nu>\frac{1}{2}n+1$, 
then $\Pic^n_{\nu}=\emptyset$. Indeed, Clifford's theorem gives that for such couples $(n,\nu)$
we have $\Pic^n_{\nu}(\overline{K})=\emptyset$ and as $\Pic^n_{\nu}$ is of finite
type over $K$ this implies $\Pic_{\nu}^n=\emptyset$. The rest is left to the reader.
\end{proof}

There is the following irreducibility result for the two-variable motivic zeta function
which contains the analogous result for Pellikaans two-variable zeta function as a special case:

\begin{theorem}

Let $K$ be a field and $X/K$ a smooth projective geometrically irreducible curve
of genus at least one. Assume $\mu:K_0(\Var_K)\rightarrow A$ is a motivic measure with
values in a field $A$. Then $P_{X,\mu}(T,u)$ is (absolutely) irreducible 
if and only if $\mu(\Pic^0_{X/K})\neq 0$.

\end{theorem}

\begin{proof}

If $\mu:=\mu(\Pic^0_{X/K})=0$ then $P_{X,\mu}(1,u)=0$ by theorem \ref{functionalequation}, d) so $P_{X,\mu}(T,u)$
is divisible by $T-1$. 

The proof of absolute irreducibility of $P_{X,\mu}(T,u)$ in case $\mu\neq 0$ is virtually the 
same as for Pellikaan's function, which can be found in \cite{N}. So we only give a short
sketch here:
We prove absolute irreducibility for $F(T,u):=T^{2g}P_{X,\mu}(T^{-1},u)$ which implies the result because $P_{X,\mu}(0,u)=1\neq 0$. The polynomial $F$ is monic in $T$ and satisfies $F(1,u)=\mu$
a non-zero {\em constant}. An easy exercise in commutative algebra shows that to finish the proof it is sufficient
to show that the leading coefficient of $F$ as a polynomial {\em in u} is irreducible in
$A[T]$. Indeed, we have more precisely:

\[
F(T,u)=(1-T)u^g+O(u^{g-1}).
\]

Using $\deg\; P_i(u)\leq\frac{1}{2}+i$ this boils down to an assertion about $P_1(u)$ and
 $P_2(u)$ for which we refer to {\em loc.cit.} where one has to replace $b_{nk}$ by
$\mu(\Pic^n_k)$. A final point is that the vanishing of certain $b_{nk}$ deduced from
Clifford's theorem in {\em loc.cit.} carries over to the $\mu(\Pic^n_k)$ as explained at the end of the proof of theorem \ref{functionalequation}.

\end{proof}

authors' adresses:\\

Francesco Baldassarri\\
Dipartimento di Matematica Pura et Applicata \\
Via G. Belzoni 7\\
35131 Padova\\
Italy\\
email: baldassa@math.unipd.it\\

Christopher Deninger\\
Mathematisches Institut der WWU M\"unster\\
Einsteinstr. 62\\
48149 M\"unster\\
Germany\\
email: deninger@math.uni-muenster.de\\

Niko Naumann\\
Mathematisches Institut der WWU M\"unster\\
Einsteinstr. 62\\
48149 M\"unster\\
Germany\\
e-mail: naumannn@math.uni-muenster.de\\

\end{document}